\documentclass{amsart}

\usepackage{amssymb}
\usepackage{amsmath}
\usepackage{amsthm}
\usepackage{diagrams}

\newtheorem{theorem}{Theorem}[section]
\newtheorem{lemma}[theorem]{Lemma}
\newtheorem{proposition}[theorem]{Proposition}
\newtheorem{corollary}[theorem]{Corollary}

\theoremstyle{definition}
\newtheorem{definition}[theorem]{Definition}

\theoremstyle{remark}
\newtheorem{remark}[theorem]{Remark}

\numberwithin{equation}{section}

\newcommand{\f}{\varphi}
\newcommand{\B}{{\mathcal B}}
\newcommand{\V}{{\mathcal V}}
\newcommand{\N}{\mathbb N}
\newcommand{\pu}{\mathbb P^1}
\newcommand{\proj}{\mathbb P}
\DeclareMathOperator{\loc}{\mathrm{Locus}}
\DeclareMathOperator{\cloc}{\mathrm{ChLocus}}
\newcommand{\ratcurves}{\textrm{Ratcurves}^n(X)}
\newcommand{\ratcurvesx}{\textrm{Ratcurves}^n(X,x)}
\newcommand{\om}{\textrm{Hom}}
\newcommand{\Aut}{\textrm{Aut}}
\newcommand{\Univ}{\textrm{Univ}}

\begin{document}

\title[Generalized Mukai conjecture for special Fano varieties]{Generalized Mukai conjecture\\ for special Fano 
varieties}

\author{Marco Andreatta}
\address{Dipartimento di Matematica,  Universit\`a di Trento, via Sommarive 14, I-38050 Povo (TN)}
\email{andreatt@science.unitn.it}

\author{Elena Chierici}
\address{Dipartimento di Matematica,  Universit\`a di Trento, via Sommarive 14, I-38050 Povo (TN)}
\email{chierici@science.unitn.it}

\author{Gianluca Occhetta}
\address{Dipartimento di Matematica,  Universit\`a di Trento, via Sommarive 14, I-38050 Povo (TN)}
\email{occhetta@science.unitn.it}



\begin{abstract}
Let $X$ be a Fano variety of dimension $n$, pseudoindex $i_X$ and
Picard number $\rho_X$. A generalization of a conjecture of Mukai says that $\rho_X(i_X-1)\le n$.
We prove that the conjecture holds if: a) $X$ has pseudoindex $i_X \ge \frac{n+3}{3}$
and either has a fiber type extremal contraction or does not have small extremal contractions
b) $X$ has dimension five.
\end{abstract}

\maketitle

\section{Introduction}

Let $X$ be a Fano variety, that is a smooth complex projective variety whose
anticanonical bundle $-K_X$ is ample. We denote with $r_X$ the {\sf index} of $X$ and with $i_X$
the  {\sf pseudoindex} of $X$, defined respectively as
$$r_X = \max \{m \in \N\ |-K_X = mL  {\mathrm{\ for \ some \ line \ bundle \ }} L  \},$$
$$i_X = \min \{m \in \N\ |-K_X \cdot C = m {\mathrm{\ for \ some \ rational \ curve \ }} C \subset X \}.$$

\medskip
In 1988, Mukai \cite{Kata} proposed the following conjecture:

\smallskip
{\bf Conjecture A.}\quad
Let $X$ be a Fano variety of dimension $n$. Then
   $$\rho_X(r_X - 1) \le n.$$

A more general conjecture (since $i_X \geq r_X$), which we will consider here,
has the following form:

\smallskip
{\bf Conjecture B.}\quad
Let $X$ be a Fano variety of dimension n. Then
   $$\rho_X(i_X-1) \le n,$$
with equality if and only if $X \simeq (\proj^{i_X-1})^{\rho_X}$.\par

\medskip
In 1990 Wi\'sniewski \cite{Wimu} proved that if $i_X > \frac{n+2}{2}$ then $\rho_X =1$;
in that paper he implicitly noticed that the statement of Conjecture B is more natural.
In 2002 Bonavero, Casagrande, Debarre and Druel \cite{BCDD} explicitely posed conjecture B
and proved it in the following situations: (a) $X$ has dimension $4$, (b) $X$ is a toric
variety of pseudoindex $i_X \geq \frac{n+3}{3}$ or of dimension $\le 7$. In this paper we prove the following

\begin{theorem}\label{main}
Let $X$ be a Fano variety of dimension $n$; then conjecture B
holds in the following cases:
   \begin{enumerate}
     \item[(a)] $i_X \ge \frac{n+3}{3}$ and $X$ has a fiber type extremal contraction;
	 \item[(b)] $i_X \ge \frac{n+3}{3}$ and $X$ has not small extremal contractions;
	 \item[(c)] $n = 5$.
   \end{enumerate}
\end{theorem}

We use the language of the Minimal Model Program, or Mori theory; 
therefore for us an extremal contraction is a map with connected fibers from $X$ onto a 
normal projective variety; such a map contracts 
all curves in an extremal face of the Kleiman-Mori cone $\overline{NE(X)} \subset N_1(X)$. 
Remember that, since $X$ is Fano, $NE(X)$ 
is contained in the  half space defined by $\{ z \in N_1(X): K_X z < 
0\}$ and so, by the Cone theorem,  $NE(X)$ is a polyhedral closed 
cone.\par
\medskip
Note that, while condition (b) is certainly a strong one, condition (a) seems very natural;
actually we do not know any example of a Fano variety which does not have fiber type
contractions.\par
\medskip
We use the typical tools for this kind of problems, in particular the existence of ``many''
rational curves on $X$ which is a fundamental property of Fano varieties as shown by Mori
\cite{Mo2}.\\
We work with families of rational curves, i.e. components of the scheme \linebreak $\ratcurves$ which
parametrizes birational morphisms $\pu \to X$ up to automorphisms of $\pu$, and families
of rational 1-cycles, i.e. components of $\textrm{Chow}(X)$, which we call Chow families;
we will denote families of rational curves by capital letters and Chow families by
calligraphic letters.\\
To a family of rational curves $V$ one can associate a Chow family $\V$, taking the closure
of the image of $V$ in $\textrm{Chow}(X)$ via the natural morphism \linebreak
$\ratcurves \to \textrm{Chow}(X)$; if $V$ is an unsplit family, i.e. if $V$ is a proper scheme,
then the two notions essentially agree and we can identify $V$ with $\V$.\par

\medskip
Fano varieties are rationally connected, i.e. through every pair of points $x,y \in X$ there
exists a rational curve; this was proved in \cite{Cam} and in \cite{KoMiMo}.
In this paper, as in \cite{AW}, we use the notion of rational connectedness
with respect to some chosen Chow families of rational curves $\V^1, \dots, \V^k$: roughly speaking,
$X$ is rc$(\V^1, \dots, \V^k)$ connected if through every pair of points $x,y \in X$ there
passes a connected 1-cycle whose components belong to the families $\V^1, \dots, \V^k$.\\
To the rc$(\V^1, \dots, \V^k)$ relation one can associate a proper fibration, called
rationally connected fibration, defined on an open set of $X$, whose fibers are
equivalence classes for the relation; this was proved again in \cite{Cam} and \cite{KoMiMo}.\par

\medskip
Using this fact we prove that if $X$ is rationally connected with respect to $k$ unsplit families 
$V^1, \dots, V^k$ then $\rho_X \le k$.\\
Then we show that if $X$ satisfies assumption (a) or (b) of theorem \ref{main} 
then $X$ is rationally connected with respect
to $k \le 3$ unsplit families, and equality holds if and only if $X = (\proj^{i_X-1})^3$.\\
The case of Fano fivefolds is more difficult: we prove that $X$ is rationally connected
with respect to a suitable number of proper families, but one of them could be
a non unsplit Chow family, so to get the result we have to bound the number of its
possible splittings.


\section{Families of rational curves}

We recall some of our basic definitions; our notation is basically consistent with the one in
\cite{Kob} to which we refer the reader.\\
Let $X$ be a normal projective variety and let $\om(\pu,X)$ be the scheme parametrizing
morphisms $f: \pu \to X$; we consider $\om_{bir}(\pu,X) \subset \om(\pu,X)$, the open subscheme
corresponding to those morphisms which are birational onto their image, and its normalization
$\om^n_{bir}(\pu,X)$; the group $\Aut(\pu)$ acts on $\om^n_{bir}(\pu,X)$ and the quotient exists.

\begin{definition}
The space $\ratcurves$ is the quotient of $\om^n_{bir}(\pu,X)$ by $\Aut(\pu)$, and the space
$\Univ(X)$ is the quotient of the product action of $\Aut(\pu)$ on
$\om^n_{bir}(\pu, X)\times\pu$.
\end{definition}

We have the following commutative diagram:
\begin{diagram}[loose,height=2em,w=4em,l>=4em]
   \om^n_{bir} (\pu, X)\times \pu  & \rTo^U & \Univ(X)\\
   \dTo & & \dTo>p\\
   \om^n_{bir} (\pu, X) & \rTo^u & \ratcurves
\end{diagram}
where $u$ and $U$ are principal $\Aut(\pu)$-bundles and $p$ is a $\pu$-bundle.\par

\begin{definition} \label{Rf}
We define a {\sf family of rational curves} to be an irreducible component
$V \subset \ratcurves$.\\
Given a rational curve $f:\pu \to X$ we will call a {\sf family of deformations}
of $f$ any irreducible component $V \subset \ratcurves$ containing $u(f)$.
\end{definition}

Given a family of rational curves, we have the following basic diagram:
\begin{diagram}[size=2em] \label{diagram}
   p^{-1}(V)=:U & \rTo^{i} & X \\
   \dTo<p & & \\
   V & &\\
\end{diagram}
where $i$ is the map induced by the evaluation $ev:\om^n_{bir} (\pu, X)\times \pu \to X$
and $p$ is a $\pu$-bundle. We define $\loc(V)$ to be the image of $U$ in $X$;
we say that $V$ is a {\sf covering family} if $i$ is dominant, i.e. if $\overline{\loc(V)}=X$.
We will denote by $\deg V$ the anticanonical degree of the family $V$, i.e. the integer
$-K_X \cdot C$ for any curve $C \in V$.\\
If we fix a point $x \in X$, everything can be repeated starting from the scheme
$\om(\pu,X; 0 \mapsto x)$ which parametrizes morphisms $f:\pu \to X$ sending $0 \in \pu$ to $x$. 
Again we obtain a commutative diagram
\begin{equation}
   \begin{diagram}[loose,height=2em,w=4em,l>=4em]
      \om^n_{bir} (\pu, X; 0 \mapsto x)\times \pu  & \rTo^U & \Univ(X,x)\\
      \dTo & & \dTo>p\\
      \om^n_{bir} (\pu, X; 0 \mapsto x) & \rTo^u & \ratcurvesx
   \end{diagram}
\end{equation}
and, given a family $V \subseteq \ratcurves$, we can consider the subscheme \linebreak
$V \cap \ratcurvesx$ parametrizing curves in $V$ passing through $x$. We usually denote
by $V_x$ a component of this subscheme.

\begin{definition} Let $V$ be a family of rational curves on $X$. Then
\begin{enumerate}
\item[(a)] $V$ is {\sf unsplit} if it is proper;
\item[(b)] $V$ is {\sf locally unsplit} if for the general $x \in \loc(V)$ every component $V_x$ is proper.
\item[(c)] $V$ is {\sf generically unsplit} if there is
    at most a finite number of curves of $V$ passing through two general points of $\loc(V)$.
\end{enumerate}
\end{definition}

\begin{remark} Note that $(a) \Rightarrow (b) \Rightarrow (c)$.
\end{remark}

\begin{proposition}\cite[IV.2.6]{Kob}\label{iowifam} Let $X$ be a smooth projective variety 
and let $V$ be a family of rational curves.\\
Assume either that $V$ is generically unsplit and $x$ is a general point in $\loc(V)$ or that
$V$ is unsplit and $x$ is any point in $\loc(V)$. Then
 \begin{enumerate}
      \item[(a)] $\dim X + \deg V \le \dim \loc(V)+\dim \loc(V_x) +1=\dim V$;
      \item[(b)] $\deg V \le \dim \loc(V_x)+1$.
   \end{enumerate}
\end{proposition}

\begin{definition}\label{CF}
We define a {\sf Chow family of rational curves} to be an irreducible component
$\V \subset \textrm{Chow}(X)$ parametrizing rational connected 1-cycles.
\end{definition}

Given a Chow family of rational curves, we have a diagram as before, coming from the universal
family over $\textrm{Chow}(X)$.
\begin{equation}\label{chowdiagram}
   \begin{diagram}[size=2em]
      {\mathcal U} & \rTo^{i} & X \\
      \dTo<p & & \\
      \V & &\\
   \end{diagram}
\end{equation}
In the diagram $i$ is the map induced by the evaluation and the fibers of $p$ are connected and have
rational components. Both $i$ and $p$ are proper (see for instance \cite[II.2.2]{Kob}).
By \cite[IV.4.10]{Kob} the family $\V$ defines a proper prerelation in the sense of
\cite[IV.4.6]{Kob} (note that schemes and morphisms appearing in that definition are those
of the normal form \cite[IV.4.4.5]{Kob}).

\begin{definition}
If $V$ is a family of rational curves we can consider the closure of the image of
$V$ in $\textrm{Chow}(X)$, and call it the {\sf Chow family associated to} $V$.
\end{definition}

\begin{remark}
If $V$ is proper, i.e. if the family is unsplit, then $V$ corresponds to the normalization
of the associated Chow family $\V$; in particular $V$ itself defines a proper prerelation.
\end{remark}


\section{Chains of rational curves}

Let $X$ be a normal proper variety, $\V^1, \dots, \V^k$ Chow families of rational curves
on $X$ and $Y$ a subset of $X$.

\begin{definition}
We denote by $\loc(\V^1, \dots, \V^k)$ the set of points $x \in X$ such that there exist
cycles $C_1, \dots, C_k$ with the following properties:
   \begin{itemize}
      \item $C_i$ belongs to the family $\V^i$;
      \item $C_i \cap C_{i+1} \neq \emptyset$;
      \item $x \in C_1 \cup \dots \cup C_k$,
   \end{itemize}
i.e. $\loc(\V^1, \dots, \V^k)$ is the set of points which belong to a connected chain of
$k$ cycles belonging \underline{respectively} to the families $\V^1, \dots, \V^k$.
\end{definition}

Note that if $\V$ is a Chow family then $\loc(\V)$ is the image of $\mathcal U$ in $X$
through $i$ in diagram \ref{chowdiagram}, so, since $\V, p$ and $i$ are proper,
$\loc(\V)$ is a closed subset of $X$.

\begin{definition}
We denote by $\loc(\V^1, \dots, \V^k)_Y$ the set of points $x \in X$ such that there exist
cycles $C_1, \dots, C_k$ with the following properties:
   \begin{itemize}
      \item $C_i$ belongs to the family $\V^i$;
      \item $C_i \cap C_{i+1} \not = \emptyset$;
      \item $C_1 \cap Y \not = \emptyset$ and $x \in C_k$,
   \end{itemize}
i.e. $\loc(\V^1, \dots, \V^k)_Y$ is the set of points that can be joined to $Y$ by a connected
chain of $k$ cycles belonging \underline{respectively} to the families $\V^1, \dots, \V^k$.
\end{definition}

Note that $\loc(\V^1, \dots, \V^k)_Y \subset \loc(\V^k)$.

\begin{remark} \label{Lr?}
If $Y$ is a closed subset, then $\loc(\V^1, \dots, \V^k)_Y$ is closed.\\
We prove the statement by induction, since we have
   $$\loc(\V^1, \dots, \V^k)_Y= \loc(\V^k)_{\loc(\V^1, \dots, \V^{k-1})_Y}.$$
With the notation of diagram \ref{chowdiagram} let $\V_Y= p(i^{-1}(Y \cap \loc(\V)))$
be the subset of $\V$ parametrizing cycles of $\V$ meeting $Y$;
$\loc(\V)_Y$ is just $i(p^{-1}(\V_Y))$, so it is closed by the properness of $i$ and $p$.
\end{remark}

\begin{definition}
We denote by $\cloc_m(\V^1, \dots, \V^k)_Y$ the set of points $x \in X$
such that there exist cycles $C_1, \dots, C_m$ with the following properties:
   \begin{itemize}
      \item $C_i$ belongs to a family $\V^j$;
      \item $C_i \cap C_{i+1} \not = \emptyset$;
      \item $C_1 \cap Y \not = \emptyset$ and $x \in C_m$,
   \end{itemize}
i.e. $\cloc_m(\V^1, \dots, \V^k)_Y$ is the set of points that can be joined to $Y$ by a
connected chain  of at most $m$ cycles belonging to the families $\V^1, \dots, \V^k$.
\end{definition}

\begin{remark} \label{chainunion}
Note that
   $$\cloc_m(\V^1,\dots,\V^k)_Y =\bigcup_{1\leq i(j)\leq k} \loc(\V^{i(1)},\dots,\V^{i(m)})_Y;$$
in particular, if $Y$ is a closed subset then $\cloc_m(\V^1, \dots, \V^k)_Y$ is closed.
\end{remark}

\begin{definition}
We define a relation of {\sf rational connectedness with respect to $\V^1, \dots, \V^k$}
on $X$ in the following way: $x$ and $y$ are in rc$(\V^1,\dots,\V^k)$ relation if there
exists a chain of rational curves in $\V^1, \dots ,\V^k$ which joins $x$ and $y$, i.e.
if $y \in \cloc_m(\V^1, \dots, \V^k)_x$ for some $m$.
\end{definition}

\begin{remark}
The rc$(\V^1,\dots,\V^k)$ relation is nothing but the set theoretic relation
$\langle {\mathcal U}_1, \dots, {\mathcal U}_k \rangle$ associated to the proper proalgebraic
relation $\textrm{Chain}({\mathcal U}_1, \dots, {\mathcal U}_k)$ in the language of
\cite[IV.4.8]{Kob}.
\end{remark}

To the rc$(\V^1,\dots,\V^k)$ relation we can associate a fibration, at least on an open subset.

\begin{theorem}\cite[IV.4.16]{Kob}\label{rcvfibration}
There exist an open subvariety $X^0 \subset X$ and a proper morphism with connected fibers
$\pi:X^0 \to Z^0$ such that
   \begin{enumerate}
      \item[(a)] the rc$(\V^1,\dots,\V^k)$ relation restricts to an equivalence relation on $X^0$;
      \item[(b)] the fibers of $\pi$ are equivalence classes for the rc$(\V^1,\dots,\V^k)$ relation;
      \item[(c)] for every $z \in Z^0$ any two points in $\pi^{-1}(z)$ can be connected by a chain
      of at most $2^{\dim X - \dim Z}-1$ cycles in $\V^1, \dots, \V^k$.
   \end{enumerate}
\end{theorem}


\section{Bounding the Picard number of $X$}

\begin{lemma}\label{numeq}
Let $Y \subset X$ be a closed subset, $\V$ a Chow family of rational curves. Then every curve
contained in $\loc(\V)_Y$ is numerically equivalent to a linear combination with rational
coefficients of a curve contained in $Y$ and irreducible components of cycles parametrized
by $\V$ which intersect $Y$.
\end{lemma}

\begin{proof}
Let $\V_Y=p(i^{-1}(Y \cap \loc(\V)))$, let ${\mathcal U}_Y =p^{-1}(V_Y)$
and consider the restriction of diagram \ref{chowdiagram}
\begin{diagram}[size=2em]
   {\mathcal U}_Y & \rTo^{i} & X \\
   \dTo<p & & \\
   \V_Y & &\\
\end{diagram}
Let $C$ be a curve in $\loc(\V)_Y$ which is not an irreducible component of a cycle
parametrized by $\V$. Then $i^{-1} (C)$ contains an irreducible curve $C'$ which is not
contained in a fiber of $p$ and dominates $C$ via $i$.
Let $B = p(C')$ and let $S$ be the surface $ p^{-1}(B)$.\\
Note that there is a curve $C'_Y$ in $S$ which dominates $B$ and such that $i(C'_Y)$ is
contained in $Y$: this is due to the fact that the image via $i$ of every fiber of $p_{|S}$
meets $Y$.\\
By \cite[II.4.19]{Kob} every curve in $S$ is algebraically equivalent to a linear combination
with rational coefficients of $C'_Y$ and of the irreducible components of fibers of $p_{|S}$
(in \cite[II.4.19]{Kob} take $X = S$, $Y= B$ and $Z = C'_Y$).\\
Thus any curve in $i(S)$, and in particular $C$, is algebraically, hence numerically,
equivalent in $i({\mathcal U}_Y)= \loc(\V)_Y$ (and hence in $X$) 
to a linear combination with rational coefficients of $i_*(C_Y)$ and of
irreducible components of cycles parametrized by $\V_Y$.\end{proof}

\begin{corollary}\label{chain}
Let $Y\subset X$ be a closed subset, $\V^1, \dots, \V^k$ Chow families of rational curves,
$m$ a positive integer.\\
Then every curve contained in $\cloc_m(\V^1, \dots, \V^k)_Y$ is numerically equivalent
to a linear combination with rational coefficients of a curve contained in $Y$ and
irreducible components of cycles in $\V^1, \dots, \V^k$.
\end{corollary}

\begin{proof}
By \ref{chainunion}, $\cloc_m(\V^1, \dots, \V^k)_Y= \bigcup_{1 \leq i(j) \leq k} \loc(\V^{i(1)},
\dots,\V^{i(m)})_Y,$ so every irreducible component of $\cloc_m(\V^1, \dots, \V^k)_Y$ is
contained in \linebreak
$\loc(\V^{i(1)},\dots,\V^{i(m)})_Y$ for some $m$-uple $(i(1), \dots, i(m))$.\\
Then we note that the corollary is true for $\loc(\V^{i(1)}, \dots, \V^{i(m)})_Y$, applying
$m$ times lemma \ref{numeq} with $Y_0=Y$ and $Y_j =\loc(\V^{i(1)}, \dots, \V^{i(j)})_Y$. \end{proof}

\begin{proposition}\label{rhobound2}
Let $\V^1,\dots, \V^k$ be Chow families of rational curves on $X$ and let $\pi:X^0 \to Z^0$
be the rc$(\V^1,\dots, \V^k)$ fibration.\\
Let $Y \subset X$ be a closed subset which dominates $Z^0$ via $\pi$; then every curve in $X$
is numerically equivalent to a linear combination with rational coefficients of a curve
contained in $Y$ and irreducible components of cycles in $\V^1, \dots, \V^k$.
\end{proposition}

\begin{proof}
By theorem \ref{rcvfibration} and the assumption, every couple of points in a general fiber
of $\pi$ can be connected by a chain of cycles belonging to $\V^1, \dots, \V^k$ of length at
most $M=2^{\dim X-\dim Z}-1$. In particular it follows that $\cloc_M(\V^1, \dots, \V^k)_Y$ is
dense in $X$ and, being closed by remark \ref{chainunion}, it coincides with $X$. Then the
claim follows from corollary \ref{chain}.\end{proof}

\begin{corollary}\label{rhobound}
Suppose that $X$ is rationally connected with respect to some Chow families
$\V^1,\dots, \V^k$; then every curve in $X$ is numerically equivalent to a linear combination with rational
coefficients of irreducible components of cycles in $\V^1, \dots, \V^k$. In particular if $X$ is rationally 
connected with respect to $k$ unsplit families, then $\rho_X \le k$.
\end{corollary}

\begin{proof} We apply proposition \ref{rhobound2} with $\pi: X \to \{*\}$ the contraction of $X$
to a point and $Y$ any point in $X$. The second part follows from the fact that any cycle parametrized by 
an unsplit family is irreducible. \end{proof}


\section{Unsplit families}

The results in the previous section can be enriched if we consider unsplit families of rational
curves instead of Chow families.

\begin{lemma} \cite[Lemma 1]{op} \label{numequns}
Let $Y \subset X$ be a closed subset, $V$ an unsplit family of rational curves. Then
$\loc(V)_Y$ is closed and every curve contained in $\loc(V)_Y$ is numerically equivalent to
a linear combination with rational coefficients
   $$\lambda C_Y + \mu C_V,$$
where $C_Y$ is a curve in $Y$, $C_V$ belongs to the family $V$ and $\lambda \ge 0$.
\end{lemma}

Note that the improvement with respect to lemma \ref{numeq} is the claim $\lambda \ge 0$.

\begin{corollary}\label{ray} Let $R=\mathbb R_+[\Gamma]$ be an extremal ray of $X$,
$V_\Gamma$ a family of deformations of a minimal extremal curve, $x$ a point in $\loc(V_\Gamma)$ and 
$V$ an unsplit family of rational curves, independent from $V_\Gamma$.\\ 
Then every curve contained in $\loc(V_\Gamma,V)_x$  is numerically equivalent to a 
linear combination with rational coefficients
   $$\lambda C_V + \mu C_\Gamma,$$
where $C_V$ is a curve in $V$, $C_\Gamma$ belongs to the family $V_\Gamma$ and $\lambda,\mu \ge 0$.
\end{corollary}

\begin{proof} By lemma \ref{numequns}, if $C$ is a curve in 
$\loc(V_\Gamma,V)_x = \loc(V)_{\loc(V_\Gamma)_x}$, then
 $$C \equiv \lambda C_\Gamma + \mu C_V,$$
with $\lambda \ge 0$ so we have only to prove that $\mu \ge 0$.\\
If $\mu <0$, then we can write $C_\Gamma \equiv \alpha C_V + \beta C$ with $\alpha, \beta \ge 0$, but
since $C_\Gamma$ is extremal, this implies that both $[C]$ and $[C_V]$ belong to $R$, a contradiction.\end{proof} 
\bigskip 

\bigskip  
One of the advantages of using unsplit families is given by the existence of good estimates
for the dimension of $\loc(V^1, \ldots, V^k)_x$:

\begin{theorem} \cite[Th\'eor\`eme 5.2]{BCDD} \label{bcdd}
Let $V^1, \ldots, V^k$ be unsplit families of rational curves on $X$. If the corresponding
classes in $N_1(X)$ are independent, then either $\loc(V^1, \ldots, V^k)_x$ is empty or it has
dimension greater or equal to $\sum \deg V^i - k$.
\end{theorem}

Using the same techniques as in the proof of theorem \ref{bcdd} we obtain the following:

\begin{lemma} \label{locy}
Let $Y \subset X$ be a closed subset and $V$ an unsplit family.  
Assume that curves contained in $Y$ are numerically independent from curves in $V$, and that 
$Y \cap \loc(V) \not= \emptyset$. Then for a general $y \in Y \cap \loc(V)$ 
\begin{enumerate}
      \item[(a)] $\dim \loc(V)_Y \ge \dim (Y \cap \loc(V)) + \dim \loc(V_y);$
      \item[(b)] $\dim \loc(V)_Y \ge \dim Y + \deg V - 1$.
\end{enumerate}	  
Moreover, if $V^1, \dots, V^k$ are numerically independent 
unsplit families such that curves contained in $Y$ are numerically independent 
from curves in $V^1, \dots, V^k$ then either $\loc(V^1, \ldots, V^k)_Y=\emptyset$
or
\begin{enumerate}
	  \item[(c)] $\dim \loc(V^1, \ldots, V^k)_Y \ge \dim Y +\sum \deg V^i -k$.
\end{enumerate}
\end{lemma}

\begin{proof}
We refer to diagram  \ref{diagram}. Since $V$ is unsplit, for a point $y$ in $Y \cap \loc(V)$
we have
   $$\dim i^{-1}(y) = \dim V_y = \dim \loc(V_y) -1.$$
So, setting $V_Y = p(i^{-1}(Y))$ and $U_Y =p^{-1}(V_Y)$, we have for general \linebreak
$y \in Y \cap \loc(V)$,
\begin{eqnarray*}
   \dim U_Y  & = & \dim (Y \cap \loc(V))+ \dim \loc(V_y)  \ge  \\
   & \ge & \dim Y + \dim \loc(V) - n +\dim \loc(V_y) \ge \\
   & \ge & \dim Y + \deg V - 1.
\end{eqnarray*}
Since $\loc(V)_Y =i(U_Y)$, (a) and (b) will follow if we prove that $i: U_Y \to X$ is generically finite.\\
To show this we take a point $x \in i(U_Y) \backslash Y$ and we suppose that
$i^{-1}(x) \cap U_Y$ contains a curve $C'$ which is not contained in any fiber of $p$; let
$B'$ be the curve $p(C') \subset V_Y$ and let $\nu:B \to B'$ be the normalization of $B'$.\\
By base change we obtain the following diagram
\begin{diagram}[height=2em,width=2.5em]
   S_B & \rTo^{j} &X\\
   \dTo<{p_B} & &\\
   B & &
\end{diagram}
Let $C_Y$ be a curve in $S_B$ which dominates $B$ and whose image via $j$ is contained in $Y$;
such a curve exists since the image via $j$ of every fiber of $p_B$ meets $Y$.
Now two cases are possible:
either $j(C_Y)$ is a point, and therefore we have a one-parameter family of curves passing
through two fixed points, contradicting the fact that $V$ is unsplit
(see for instance \cite[IV.2.3]{Kob})
or $j(C_Y)$ is a curve in $Y \cap \loc(V_y)$, so a curve in $Y$ is numerically proportional to
a curve parametrized by $V$, against the assumptions.\\
To show (c) it is enough to recall that, as already observed in remark \ref{Lr?}, we have 
$\loc(\V^1, \dots, \V^k)_Y= \loc(\V^k)_{\loc(\V^1, \dots, \V^{k-1})_Y}.$\end{proof}

\begin{remark}\label{noncovering} 
If in the previous theorem $V^1$ is not a covering family and \linebreak
$\loc(V^1, \ldots, V^k)_x$ is nonempty, then 
  $$\dim \loc(V^1, \ldots, V^k)_x \ge \sum \deg V^i - k +1.$$
In fact $\loc(V^1, \ldots, V^k)_x= \loc(V^2, \ldots, V^k)_{\loc(V^1_x)}$, and we can apply part (c) 
of lemma \ref{locy}, recalling that $\dim \loc(V^1_x) = \deg V^1 -1$
implies that $V^1$ is covering (see proposition \ref{iowifam}).  
\end{remark}


\section{Rational curves on Fano varieties}

The geometry of Fano varieties is strongly related to the properties of families of rational
curves of low degree. The first result in this direction is a fundamental theorem, due to Mori:

\begin{theorem}\cite{Mo2}
Through every point of a Fano variety $X$ there exists a rational curve of anticanonical degree
$\leq \dim X +1$.
\end{theorem}

\begin{remark}\label{mdf}
The families $\{V^i \subset \ratcurves \}$ containing rational curves with degree $\leq n+1$
are only a finite number, so for at least one index $i$ we have that $\loc(V^i)=X$. Among these
families we choose one with minimal anticanonical degree, and call it a {\sf minimal dominating
family}. Note that every such family is locally unsplit.
\end{remark}

A relative version of Mori's theorem is the following

\begin{theorem}\cite[Theorem 2.1]{KoMiMo}
Let $X$ be a Fano manifold. Suppose that there exist a nonempty open subset $X^0$ of $X$,
a smooth quasiprojective variety of positive dimension $Z^0$ and a proper surjective morphism
$\pi:X^0 \to Z^0$. Let $z$ be a general point on $Z^0$. Then there exists a rational curve
$C$ on $X$ satisfying
   \begin{enumerate}
      \item[(a)] $C \cap \pi^{-1}(z) \neq \emptyset$;
      \item[(b)] $C$ is not contained in $\pi^{-1}(z)$;
      \item[(c)] $-K_X \cdot C \leq n+1$.
   \end{enumerate}
\end{theorem}

\begin{remark}\label{mhdf}
The families $\{V^i \subset \ratcurves \}$ containing the horizontal curves with degree
$\leq n+1$ are only a finite number, so for at least one index $i$ we have that $\loc(V^i)$
dominates $Z^0$. Among these families we choose one with minimal anticanonical degree,
and call it a {\sf minimal horizontal dominating family} for $\pi$.
\end{remark}

A typical situation where these morphisms arise is the construction of rationally connected
fibrations associated to families of rational curves, as we have explained in section 2, or more
generally to a finite number of proper connected prerelations, as done in \cite[IV.4.16]{Kob}.

\begin{lemma} \label{horizontal}
Let $X$ be a Fano variety, and let $\pi:X \rDashto Z$ be the rationally connected fibration
associated to $m$ proper connected prerelations on $X$; let $V$ be a minimal horizontal
dominating family for $\pi$. Then
   \begin{itemize}
      \item[(a)] curves parametrized by $V$ are numerically independent from curves contracted
      by $\pi$;
      \item[(b)] $V$ is locally unsplit;
      \item[(c)] if $x$ is a general point in $\loc(V)$ and $F$ is the fiber containing $x$, then
         $$\dim (F \cap \loc(V_x))=0.$$
   \end{itemize}
\end{lemma}

\begin{proof}
(a) Since $X$ is normal and $Z$ is proper, the indeterminacy locus $E$ of $\pi$ in $X$ has codimension
$\ge 2$ \cite[1.39]{De}. Pull back an ample divisor from $Z$ and
observe that it is zero on curves contracted by $\pi$.
On the other hand it intersects nontrivially curves
which are not contracted by $\pi$ and are not contained in $E$, like curves of
$V$, since $V$ is dominant.\\
(b) If for the general $x \in \loc(V)$ a curve in $V_x$ degenerates into a reducible cycle, then
at least one component of this cycle is horizontal, otherwise curves in $V$ would be numerically
equivalent to curves in the fibers. But this contradicts the minimality of $V$ among horizontal
dominating families.\\
(c) From lemma \ref{numeq} we know that any curve in $\loc(V_x)$ is numerically proportional to
$V$, while proposition \ref{rhobound2} applied to $F$ implies that all curves in $F$ can be
written as linear combinations of curves contracted by $\pi$. \end{proof}

\begin{corollary}\label{dimmhdf}
Let $X$ be a Fano variety, and let $\pi:X \rDashto Z$ be the rationally connected fibration
associated to $m$ proper connected prerelations on $X$; let $V$ be a minimal horizontal
dominating family for $\pi$. Then
   $$\deg V \le \dim Z + 1.$$
\end{corollary}

\begin{proof} It follows from lemma \ref{horizontal} and the fact that 
$\dim \loc(V_x) \ge \deg V -1$. \end{proof}


\section{Special Fano varieties of high pseudoindex}

In this section we will prove Theorem \ref{main}, (a) and (b).\\
First of all we will show that Conjecture B is true for a Fano variety $X$ of pseudoindex 
$i_X \ge \frac{\dim X+3}{3}$ which has a covering unsplit family of rational curves, then
we will prove that this is the case if $X$ is as in (a) or (b).\par

\begin{proposition}\label{unsplitB} Let $X$ be a Fano variety of dimension $n$ and pseudoindex \linebreak
$i_X \ge \frac{n+3}{3}$; if there exists a family $V$ of rational curves which is unsplit 
and covering then Conjecture B is true for $X$.
\end{proposition}

\begin{proof} Consider the rc$V$ fibration $\pi:X^0 \to Z^0$: if $\dim Z^0 =0$ then $\rho_X =1$ by
corollary \ref{rhobound} and we conclude. Otherwise take a minimal horizontal dominating family
$V'$; from lemma \ref{horizontal} we know that $V'_x$ is unsplit for general $x \in \loc(V')$.
Then applying lemma \ref{locy} (b) with $Y = \loc(V'_x)$ we obtain
\begin{eqnarray*}
   n  & \ge & \dim \loc(V)_{\loc(V'_x)} \ge \\
   & \ge & \dim \loc(V'_x) + \deg V - 1 \ge \\
   & \ge & \deg V' + \deg V - 2
\end{eqnarray*}
so $\deg V' \le 2i_X - 1$ and therefore $V'$ is unsplit.\\
Take the rc$(V,V')$ fibration $\pi': X' \to Z'$: if $\dim Z' =0$ then from corollary
\ref{rhobound} we have $\rho_X = 2$ and we conclude, otherwise take a minimal dominating
family $V''$ with respect to $\pi'$.\\
For general $x \in \loc(V'')$, denote by $F$ the fiber of $\pi'$ containing $x$: then $F$ is an
equivalence class with respect to the rc$(V, V')$ relation, so $F \supseteq \loc(V, V')_y$ for
some $y$; then theorem \ref{bcdd} implies
   $$\dim F \ge \deg V + \deg V' - 2 \ge 2i_X - 2.$$
By lemma \ref{horizontal} we have $\dim (\loc(V''_x) \cap F) =0$, so
   $$n \geq \dim F + \dim \loc(V''_x) \geq 2i_X - 2 + \deg V'' -1,$$
that is
   $$\deg V'' \leq  n+3-2i_X \leq i_X.$$
This is impossible unless $\deg V=\deg V'=\deg V''=i_X$ and $\dim \loc(V_x)=\dim \loc(V'_x)=
\dim \loc(V''_x)= i_X -1$. Proposition \ref{iowifam} implies that all these families are
covering, so we can apply \cite[Theorem 1]{op} to obtain that $X \simeq (\proj^{i_X-1})^3$.\end{proof}

\begin{theorem} Let $X$ be a Fano variety of dimension $n$ and pseudoindex $i_X \ge \frac{n+3}{3}$.
If $X$ has a fiber type extremal contraction or has not small contractions then there
exists a covering unsplit family $V$ of rational curves.
\end{theorem}

\begin{proof} First of all suppose that there exists a fiber type contraction $\f: X \to W$;
let $V_\f$ be a minimal horizontal dominating family for $\f$;
from corollary \ref{dimmhdf} we know that $\deg V_\f \le \dim W +1$.
Let $F$ be a general fiber of $\f$; we have that
$$\dim F \le \dim X- \deg V_\f +1 \le 2i_X-2.$$
By adjunction we have $K_F=(K_X)_F$, so $F$ is a Fano variety;
in particular there exists a minimal dominating family $V_F$ of degree $\le \dim F+1 \le 2i_X-1$.\\
This means that through a general point of $X$ there passes a curve of degree \linebreak
$\le 2i_X -1$, and since the families of rational curves with bounded degree are a finite number, one of them 
must be covering; the bound on the degree implies that this family is also unsplit.\par

\medskip
Suppose now that all the extremal contractions of $X$ are divisorial and, by contradiction, that
there does not exist any unsplit covering family of rational curves.\\
Let $V$ be a minimal dominating family of rational curves; since
we are assuming that $V$ is not unsplit we have $\deg V \ge 2i_X$.\\ 
Consider the Chow family $\V$ associated to $V$:
since $\deg V \leq n+1 < 3i_X$, reducible cycles in $\V$ split into exactly two irreducible
components. To each one of them we associate the corresponding irreducible component of
$\ratcurves$, which is an unsplit family.\\
We denote by $\mathcal B$ the finite set of pairs of families $(W^i, W^{i+1})$ satisfying:
\begin{itemize}
   \item $[W^i]$ is numerically independent from $[W^{i+1}]$;
   \item $[W^i] + [W^{i+1}] = [V]$;
   \item $W^i$ and $W^{i+1}$ contain irreducible components of cycles of $\V$.
\end{itemize}
Consider now the rc$\V$ fibration $\pi: X^0 \to Z^0$.\par

\medskip
{\bf Claim.} \quad $\dim Z^0=0$.\par

\bigskip
Suppose by contradiction that $Z^0$ has positive dimension, and take $V'$ a minimal horizontal
dominating family for $\pi$; we know from lemma \ref{horizontal} (c) that for a general fiber $F$
we have
   $$\dim \loc(V'_x) + \dim F \leq n,$$
which implies
   $$\deg V' \leq n+1- \dim F \leq n - 2i_X +2 <i_X,$$
a contradiction which proves the claim.\par

\bigskip
As a corollary we obtain that $N_1(X)$ is generated as a vector space by the numerical classes
of the irreducible components of cycles in $\V$ (proposition \ref{rhobound2}).\\
Note that if $[V]$ is extremal in $NE(X)$, then all the irreducible components of cycles in
$\V$ are numerically proportional to $[V]$ and $\rho_X =1$, so we can assume that $[V]$
is not extremal.\\ 
Take now $R_1=\mathbb R_+[C_1]$  to be a divisorial extremal ray of $X$, let
$E_1$ be its exceptional locus and $V^1$ an unsplit family of
deformations of a minimal extremal rational curve $C_1$.\par

\smallskip
First of all we claim that $E_1 \cdot V=0$; otherwise for a general $x \in X$ the set
$\loc(V^1)_{\loc(V_x)}$ would be nonempty, so by lemma \ref{locy} and proposition
\ref{iowifam}
   $$\dim \loc(V^1)_{\loc(V_x)} \ge \dim \loc(V_x) + \deg V^1-1 \ge 3i_X-2 > \dim X.$$
In particular we find a pair $(W^1,W^2) \in \B$ such that $E_1 \cdot W^1 <0$ and
$E_1 \cdot W^2 >0$.\par
\medskip

By corollary \ref{ray}, if $[W^1] \not = [V^1]$ then the class of every curve in $\loc(V^1,W^1)_x$ 
can be written as a linear combination with positive coefficients of $[V^1]$ and $[W^1]$, 
so, for $x \in \loc(W^1) \cap \loc(W^2)$,
   $$\dim (\loc(W^2_x) \cap \loc(V^1,W^1)_x) = 0;$$
on the other hand, if $[W^1] \neq [V^1]$, by remark \ref{noncovering} we have
   $$\dim \loc(V^1,W^1)_x \ge 2i_X -1,$$ 
and therefore
   $$\dim (\loc(W^2_x) \cap \loc(V^1,W^1)_x) \ge 3i_X -2 -n > 0.$$
We thus get a contradiction, unless $[W^1]=[V^1]$.\par

\medskip
Note that this argument also shows that for all $i \not= 1,2$ we have $E_1 \cdot W^i=0$.\par
\medskip
Since $X$ is Fano and $E_1$ is effective there exists an extremal ray $R_2$ on which
$E_1$ is positive (this is due to the fact that every effective curve on a Fano manifold can be written
as a linear combination with positive coefficients of extremal curves: see
\cite[Lemma 2]{BCW}); let $E_2$ be the exceptional locus of $R_2$.\\
We repeat the same argument and we find a pair $(W^3,W^4)$ such that
$[V^2]=[W^3]$ and $E_2 \cdot W^4 > 0$.\\
If the plane $\Pi_1$ spanned in $N_1(X)$ by the classes $[V]$ and $[V^1]$ is different
from the plane $\Pi_2$ spanned by $[V]$ and $[V^2]$, then $[V^1]$, $[V^2]$ and $[W^4]$
are independent, and $\loc(W^4,V^2,V^1)_x$ is nonempty for every $x \in \loc(W^4)$.
By remark \ref{noncovering} we get
$$\dim \loc(W^4,V^2,V^1)_x \ge 3i_X -2> n,$$
a contradiction.\\
So we suppose that $\Pi_1=\Pi_2:= \Pi$ and we choose a basis of $N_1(X)$ 
formed by $[V^1],[V]$ and by classes $[W^i]$ not contained
in $\Pi$.\\
Since the divisors $E_1$ and $E_2$ are zero on all the elements of the basis but $[V^1]$,
they are proportional in $N^1(X)$; but $E_1 \cdot V^1 < 0$ and $E_2 \cdot V^1 >0$,
so $E_1 = -k E_2$ with $k>0$.
One can now compute the intersection number of $E_1$ and $E_2$ with any curve which meets $E_1 \cup E_2$
without being contained in it, and this leads to a contradiction.\end{proof}


\section{Fano fivefolds with a covering unsplit family}

This section and the following one are devoted to the proof
of Theorem \ref{main} (c).\\
Let $X$ be a Fano variety of dimension 5 and let $V \subseteq \ratcurves$ be a minimal dominating family; 
by remark \ref{mdf} we have that $\deg V \le 6$ and $V_x$ is unsplit for a general $x \in X$. \\
If $\deg V=6$ then $X= \loc(V_x)$ and $\rho_X=1$ by lemma \ref{numeq}, therefore
we can assume $\deg V \le 5$.\par
\medskip
First of all we note that if $i_X \ge 3$, then $V$ is unsplit; moreover
in this case we can apply proposition \ref{unsplitB} and obtain the result, so from now on
we assume that ${\bf i_X=2}$ (and we thus have to prove that $\rho_X \le 5$).\par
\medskip
We divide the proof into two main cases:
in this section we will deal with the case in which {\bf $V$ is unsplit}, while in the next one
we will assume that {\bf $V$ is not unsplit}.\par
\medskip
Consider the rc$V$ fibration $\pi:X^0 \to Z^0$: if $\dim Z^0 =0$ then $\rho_X =1$ by corollary \ref{rhobound} 
and we conclude; otherwise take a minimal horizontal dominating family $V'$.\par
\medskip
{\bf Case 1. }\quad Any minimal horizontal dominating family $V'$ is not unsplit.\par 
\medskip
Note that in this situation $\deg V' \ge 4$, so $\dim \loc(V'_x)\ge 3$; in particular, since $V'$ is horizontal
and dominates $Z^0$, we have also $\dim Z^0 \ge 3$.\par
\smallskip
If $\dim Z^0 =3$ take a general point $x \in \loc(V')$, so that $V'_x$ is unsplit. Note that $Y= \loc(V'_x)$
dominates $Z^0$, so we can apply proposition \ref{rhobound2} to get $\rho_X=2$.\par
\smallskip
If $\dim Z^0=4$ consider the rc$(V, \V')$ fibration $\pi':X' \to Z'$.\par

\newpage
{\bf Claim. }\quad $\dim Z' =0$.\par
Assume that this is not the case and denote by $F'$ a general fiber of $\pi'$. 
Then there exists a minimal horizontal
dominating family $V''$ satisfying
   $$\begin{array}{rcl}
        0 = \dim(F' \cap \loc(V''_x)) & \ge & \dim F' + \dim \loc(V''_x) -5  \\
        & \ge & 4 + \dim \loc(V''_x) - 5 \\
        & \ge & \deg V'' -2
     \end{array}$$
for every $x \in F' \cap \loc(V'')$.
Thus $\deg V''=2$ and $\dim \loc(V''_x) =1 $, so by proposition \ref{iowifam} $V''$ is covering. Since $V''$
is horizontal also with respect to the fibration $\pi$ this contradicts the minimality of $V'$, thus the
claim is proved. \par

\medskip
From corollary \ref{dimmhdf} it follows that $\deg V' \le 5$, so every reducible cycle in $\V'$
splits into exactly two irreducible components; moreover the family of deformations of each component
is unsplit and non covering because of the minimality of $V'$.\\
Consider the pairs $(W^i,W^{i+1})$ of unsplit families satisfying
   \begin{itemize}
   \item $[W^i]+[W^{i+1}] = [V']$,
   \item $W^i$ and $W^{i+1}$ contain irreducible components of a cycle in $\V'$,
   \end{itemize}
and let $\B$ be the set of these pairs.\par
\medskip
If the numerical class of every pair in $\B$ lies in the plane 
$\Pi \subseteq N_1(X)$ spanned by $[V]$ and $[V']$ then, by corollary \ref{rhobound} we have that
$\rho_X = 2$ and we are done.\\
Assume therefore by contradiction that
there exists a pair $(W^1,W^2) \in \B$ whose classes don't lie in $\Pi$, call $\Pi'$ 
the plane spanned by $[W^1]$ and $[W^2]$ and set 
   $$\B_{\Pi,\Pi'} = \{(W^i,W^{i+1}) \in \B\ |\ [W^i],\  [W^{i+1}] \in <\Pi,\Pi'>~
   \textrm{and}~[W^i],[W^{i+1}] \not = [\lambda V]\}.$$
For every $(W^i, W^{i+1}) \in \B_{\Pi,\Pi'}$, for every cycle $C_i+C_{i+1} \in W^i + W^{i+1}$
and for every point $x \in C_i$ we consider $\loc(W^i,V,W^{i+1})_x$: by remark \ref{noncovering},
we have $\dim \loc(W^i,V,W^{i+1})_x \ge 4$; since $W^{i+1}$ is not covering  
every irreducible component of $\loc(W^i,V,W^{i+1})_x$ is an effective divisor on $X$,
which is contained in $\loc(W^{i+1})$. Since $W^{i+1}$ does not dominate $Z^0$, 
the intersection of any of these divisors with $V$ is zero. \\
We claim that the intersection of any of these divisors with $V'$ is also zero.\\
In fact, if $D= \loc(W^i,V,W^{i+1})_x$ is such that $D.V' >0$, then every curve in $V'$ 
intersects $\loc(W^{i+1})$. Since $V$ is covering we have 
   $$\loc(V)_{\loc(V'_x)} \supseteq \loc(V'_x),$$ 
so
   $$\loc(V,W^{i+1})_{\loc(V'_x)} \neq \emptyset;$$
we apply lemma \ref{locy} (c) and we obtain that $\dim \loc(V,W^{i+1})_{\loc(V'_x)} = 5,$
which implies that $W^{i+1}$ is covering, a contradiction.\\
Obviously we can repeat the same argument with $\loc(W^{i+1},V,W^i)_x$ for every $x \in C_{i+1}$, and we obtain
effective divisors which are contained in $\loc(W^i)$ and whose intersection with $V$ and $V'$ is zero.\par
\medskip
Call $T$ the union of all these divisors. Now take a point $y \in X \setminus T$; since $X$ is rc$(V,\V')$ 
connected, 
$y$ can be joined to $T$ by a chain of curves in $V$ and cycles in $\V'$.
In particular there exists a cycle $\Gamma$ either in $V$ or in $\V'$ which intersects $T$ but is not contained 
in it, and since every component of $T$ has intersection zero with $V$ and $V'$, it must be of the form
$C_3 + C_4$, with $(W^3,W^4) \in \B$ and $[W^3],[W^4] \not \in <\Pi,\Pi'>$.\\
So, up to exchange $W^3$ and $W^4$, there exists a component $D$ of $T$ such that \linebreak
$D \cdot W^3>0$; then $\loc(W^3)_D$ is nonempty and, by lemma \ref{locy} (b),
   $$\dim \loc(W^3)_D \ge \dim D + \deg W^3 - 1 \ge 5$$
and $W^3$ is covering, a contradiction.\par

\bigskip
{\bf Case 2.} \quad One minimal horizontal dominating family $V'$ is unsplit.\par 
\medskip
Consider the rc$(V,V')$ fibration $\pi':X' \to Z'$; if $Z'$ is a point then $\rho_X=2$ and
we conclude, otherwise take a minimal horizontal dominating family $V''$.\par
\smallskip
If $V''$ is not unsplit then $\deg V'' \ge 4$, so $\dim \loc(V''_x)\ge 3$; moreover,
since $\dim Z' \le 3$, $\loc(V''_x)$ dominates $Z'$.\\
Take a general point $x \in \loc(V'')$, so that $V''_x$ is unsplit and
apply proposition \ref{rhobound2} with $V,V'$ and $Y= \loc(V''_x)$to obtain $\rho_X=3$.\par
\smallskip
If $V''$ is unsplit we can take the rc$(V,V',V'')$ fibration $\pi'':X'' \to Z''$: either $Z''$
is a point or every minimal horizontal dominating family is unsplit.
We consider the new fibration and we repeat the same argument. Finally we find at most five
independent unsplit families on $X$ such that $X$ is rationally connected with respect to them,
so $\rho_X \le 5$ by corollary \ref{rhobound}.\\
If there are exactly five independent families,
then they must be covering and of degree $2$ and from \cite{op} we conclude that
$X \simeq (\pu)^5$.


\section{Fano fivefolds without a covering unsplit family}

We assume now that every minimal dominating family $V$ of $X$ is not unsplit, which implies that
$\deg V \ge 4$.\\
By the discussion at the beginning of the previous section we can also assume
that $i_X=2$ and that $\deg V \le 5$, so every reducible cycle in the associated Chow family $\V$
splits into exactly two irreducible components; moreover any family of deformations of each component
is unsplit and non covering because of the minimality of $V$.\\

Consider the pairs $(W^i,W^{i+1})$ of unsplit families satisfying
   \begin{itemize}
   \item $[W^i]+[W^{i+1}] = [V]$,
   \item $W^i$ and $W^{i+1}$ contain irreducible components of a cycle in $\V$,
   \end{itemize}
and let $\B$ be the set of these pairs. \par
\medskip
{\bf Claim.} \quad If $\deg V=5$ then $\rho_X = 1$.\par
\medskip
Assume by contradiction that $\deg V=5$ and $\rho_X \ge 2$.\\
Suppose that all the irreducible components of cycles in $\V$ are numerically
proportional to $V$, and consider the rc$\V$ fibration 
$\pi:X^0 \to Z^0$.\\
Now, either $Z^0$ is a point and in our assumptions $\rho_X=1$ by corollary \ref{rhobound}, or 
there exists a minimal horizontal dominating family $V'$; then for a general $x \in \loc(V')$, 
if $F$ is the fiber through $x$, we know from lemma
\ref{horizontal} that 
   $$\dim \loc(V'_x) + \dim F \le 5,$$
and since $\dim F \ge \deg V - 1 = 4$ we have $\deg V' - 1 \le \dim \loc(V'_x) \le 1$, forcing $\deg V'=2$ 
and $\dim \loc(V'_x)=1$; hence by proposition \ref{iowifam} $V'$ is covering, against the assumptions.\\
So there exists a pair $(W^1,W^2) \in \B$ such that $[W^1] \neq [\alpha V]$.\par
\smallskip
Let $D$ be an irreducible component of $\loc(V_x)$
for a general $x \in X$; since $V$ is locally unsplit we have $N_1(D)=<[V]>$.\\
By proposition \ref{iowifam}, $\dim D \ge \deg V -1 \ge 4$; as we are assuming 
$\rho_X \ge 2$ it cannot be $D=X$, so $D$ is an effective divisor.\\
If $D.V=0$ then $D$ would be negative on at least a family $W^i$ and so it would contain curves in $W^i$, 
contradicting the fact that $N_1(D)=<[V]>$.\\
If else $D.V >0$, then  either $D.W^1 > 0$ or $D.W^2 > 0$; but in 
this case either $\loc(W^1_x) \cap D$ or $\loc(W^2_x) \cap D$ 
would be nonempty.
Since $W^i$ is not covering we have $\dim \loc(W^i_x) \ge 2$, therefore 
$\dim (\loc(W^1_x) \cap D) \ge 1$, against the fact that $N_1(\loc(W^i_x))=<[W^i]>$.
So the claim is proved and we can assume from now on that $\deg V=4$.\par

\bigskip
Consider the rc$\V$ fibration $\pi: X^0 \to Z^0$.\par

\bigskip
{\bf Case 1} \quad $\dim Z_0 > 0$.\par

\bigskip
In this case we actually prove that $\rho_X = 2$.\\
Choose $V'$ to be a minimal horizontal dominating family for $\pi$; again we know that
   $$\dim \loc(V'_x) + \dim F \leq 5,$$
but in this case $\dim F \ge 3$ so $\deg V' - 1\le \dim \loc(V'_x) \le 2$.\\
On the other hand $\dim \loc (V'_x) \ge \deg V' \ge 2$, since otherwise 
$V'$ would  be covering and of degree $2$ by \ref{iowifam}, contradicting the minimality of $V$.\\
It follows that $\dim F=3$, $\dim \loc (V'_x)=2$ and $\deg V' = 2$, so $V'$ is unsplit 
and $\dim \loc(V')=4$.\\
Moreover, since $\loc(V'_x)$ meets the general fiber of $\pi$, then $X$ is rc$(\V,V')$ connected.\par

\medskip
Let $\Pi$ be the plane spanned by $[V]$ and $[V']$ and let 
$$\B_\Pi = \{(W^i,W^{i+1}) \in \B\ |\ [W^i]\  \mathrm{and}\  [W^{i+1}] \in \Pi\}.$$ 
If $\B_\Pi = \B$ then we have $\rho_X =2$ by corollary \ref{rhobound}.\\
Suppose that this is not the case and let $D'$ be an irreducible component of $\loc(V')$.\\
Since $D'$ does not contain the general fiber $F$ of $\pi$ and the general $F$ coincides with 
$\loc(V_x)$ for some $x$, there exists a curve of $V$ meeting $D'$ but not entirely contained in it; 
therefore $D'.V>0$.\\
Let $V'_{D'}$ be the closed subfamily of $V'$ such that $\loc(V'_{D'})=D'$; by lemma \ref{locy} (a)
$\dim \loc(V'_{D'})_{D' \cap F}\ge 4$ i.e. $\loc(V'_{D'})_{D' \cap F} =D'$.\\
Since $N_1(D' \cap F)=<[V]>$, lemma \ref{numeq} implies that $N_1(D')=<[V],[V']>$.\\ 
Let $(W^1,W^2)$ be a pair in $\B \setminus \B_\Pi$; since $D'.V >0$ 
either $D'.W^1>0$ or \linebreak
$D'.W^2>0$, so we can assume $\loc(W^1)_{D'} \neq \emptyset$; but this implies
by lemma \ref{locy} that \linebreak $\dim \loc(W^1)_{D'} \ge 5$, and therefore that 
$W^1$ is covering, a contradiction.\par

\bigskip
{\bf Case 2 } \quad $\dim Z_0= 0~~$ i.e. $X$ is rc$\V$ connected. \par 

\bigskip
In this case by corollary \ref{rhobound} $N_1(X)$ is generated as a vector space by the
numerical classes of the irreducible components of cycles in $\V$.\\
We want to show that $\rho_X \le 3$, so by contradiction we assume that there exist
three pairs $(W^1,W^2), (W^3,W^4)$ and $(W^5,W^6)$ in $\B$ whose classes generate
a four dimensional vector space inside $N_1(X)$.\\
Let $\Pi \subset N_1(X)$ be the plane generated by $[W^1]$ and $[W^2]$,
and let
$$\B_\Pi = \{(W^i,W^{i+1}) \in \B\ |\ [W^i]\  \mathrm{and}\  [W^{i+1}] \in \Pi\}.$$ 

For every pair $(W^i, W^{i+1}) \in \B_\Pi$ let $\{D^i_k\}$ be the components 
of $\loc(W^i)$ which intersect $\loc(W^{i+1})$ and let
$\{D^{i+1}_j\}$ be the components of $\loc(W^{i+1})$ which intersect $\loc(W^i)$.\\
Let us note that, by proposition \ref{iowifam}, every component of $\loc(W^i)$ has
dimension greater than three, so, since the families $W^i$ are not covering, we have
$\dim D^i_k =3$ or $4$.\par

\bigskip
{\bf Case 2a} \quad For every $i$ and every $k$ there exists $j$ such that $D^i_k = D^{i+1}_j$ 
and viceversa.\par

\bigskip
If $\dim D^i_k = 3$ then by proposition \ref{iowifam} $\dim \loc(W^i_x) = 3$ and 
$D^i_k$ is a component of $\loc(W^i_x)$ for some $x$, 
so $N_1(D_k^i) = <[W^i]>$; but since $D^i_k=D^{i+1}_j$ for some $j$, we have also $N_1(D^i_k)=<[W^{i+1}]>$,
a contradiction.\\
So we can assume that $D^i_k$ is a divisor for every $k$; 
moreover $D^i_k$ is a component of $\loc(W^i)_{\loc(W^{i+1}_x)}$, 
and so $N_1(D^i_k)=<[W^i],[W^{i+1}]>$.\par

\medskip
Let us consider the intersection number of one of these divisors, say $D^1_1=:D$, with the family $V$;
if $D \cdot V >0$ then, up to exchange $W^3$ and $W^4$, we have $D \cdot W^3 >0$.\\
By lemma \ref{locy}, since $\loc(W^3)_D$ is nonempty, we have $\dim \loc(W^3)_D =5$, 
a contradiction since $W^3$ is not covering.\par

\medskip
Therefore $D \cdot V=0$, hence $D \cdot W^i <0$ for some $i$; since $N_1(D)$ is generated
by the classes of $W^1$ and $W^2$, the class of $W^i$ must belong to the plane $\Pi$.\\
In particular for every pair $(W^i,W^{i+1}) \in \B_{\Pi}$ we have that $(D \cdot W^i)(D \cdot W^{i+1}) <0$, 
yielding that
$$D= \loc(W^i)=\loc(W^{i+1}).$$

\medskip
Let now $x$ be a point outside $D$ and let $z$ be a point of $D$; since $X$ is 
rc$\V$ connected there exists a chain of cycles in $\V$ which
connects $x$ and $z$; let $\Gamma$ be the first irreducible component of one of these chains which meets $D$.\\
Since $D \cdot V=0$ then $\Gamma$ cannot belong to $V$ or to a family which is proportional to $V$. 
Moreover, since $\Gamma \not \subset D$ then $\Gamma$ does not belong to a family whose class 
is contained in the plane $\Pi$.\\
Therefore $\Gamma$ belongs to a family $W^i$ whose class is not in $\Pi$; we can thus apply
lemma \ref{locy} and obtain $\dim \loc(W^i)_D =5,$
a contradiction, since $W^i$ is not covering. We have proved that case 2a cannot occur.\par

\bigskip
We can therefore assume, up to rename the pairs in $\B_\Pi$, that there exist meeting components 
$D_1$ and $D_2$ of $\loc(W^1)$ and $\loc(W^2)$ such that $D_1 \not = D_2$.\par

\bigskip
{\bf Case 2b} \quad $\dim D_1 = \dim D_2=4$.\par

\bigskip
We claim that we cannot have
$$D_1 \cdot V = D_2 \cdot V = 0.$$   
In fact, if $D_1 \cdot V=0$, then for at least a family $W^i$ we have $D_1 \cdot W^i <0$.\\ 
If $i \not = 1,2$ then $D_1 = \loc(W^i) = \loc(W^i)_{\loc(W^1_x)}$ and 
$N_1(D_1)=<[W^1],[W^i]>$, hence, by lemma \ref{locy}, $\dim \loc(W^2)_{D_1}=5$, a contradiction
since $W^2$ is not covering.\\
If $D_1 \cdot W^2 <0$, then $D_2 \subseteq D_1$, against our assumptions,
so we have $D_1 \cdot W^1 <0$ (and, in the same way $D_2 \cdot W^2<0$).
It follows that $D_1 = \loc(W^1) $ and $D_2=\loc(W^2)$; moreover the locus of 
every family of a pair belonging to $B_\Pi$ is contained either in $D_1$ or in $D_2$.\par

\medskip
Let $T= D_1 \cup D_2$, let $z \in T$ and let $x$ be a general point of $X$. Since $X$ is
rc$\V$ connected there exists a chain of cycles 
in $\V$ connecting $x$ and $z$; let $\Gamma$ be the first irreducible component 
which meets $T$.\\
The curve $\Gamma$ cannot be numerically proportional to $V$, since $D_1 \cdot V= D_2 \cdot V = 0$, 
and its class cannot lie in the plane $\Pi$, so $\Gamma$ belongs to an unsplit family
$W^i$ which is independent from $W^1$ and $W^2$; so either $D_1 \cdot W^i >0$ or $D_2 \cdot W^i>0$, which 
implies that either $D_1 \cdot W^{i+1}$ or $D_2 \cdot W^{i+1}$ 
is negative, a contradiction which proves the claim.\par

\medskip
Therefore we can assume that $D_1 \cdot V > 0$; up to exchange $W^3$ and $W^4$ 
we can also assume that $D_1 \cdot W^3 > 0$.\\
Let $x$ be a point on a curve in $W^4$: then $H = \loc(W^4, W^3, W^1)_x$ is nonempty and has dimension four
by remark \ref{noncovering}.\\
If $H \cdot V > 0$ then up to exchange $W^5$ and $W^6$ we can assume that $H \cdot W^5>0$, and so
by lemma \ref{locy} $\dim \loc(W^5)_H = 5,$ a contradiction.\\
If $H \cdot V = 0$, for some pair $(W^i,W^{i+1})$ we have $H. \cdot W^i < 0$ and $H \cdot W^{i+1}>0$.\\
It cannot be $i=1$, since in this case $H = \loc(W^1)= D_1$, but we are assuming that $D_1 \cdot V > 0$;
therefore $H \cdot W^i < 0$ for some $i$ such that $W^i$ is independent from $W^1$.
Let $W^1_H$ be the closed subfamily of $W^1$ whose locus is $H$;
then $H = \loc(W_H^1)_{\loc(W^i_x)}$ and so $N_1(H)=<[W^1],[W^i]>$.\\
By construction, $H \cap \loc(W^3)$ is nonempty, so either $i=3$, $H$ contains $\loc(W^3)$ 
and $\dim \loc(W^4)_H = 5$, a contradiction, or $i \neq 3$ and $\dim \loc(W^3)_H=5$, again a contradiction.
So case 2b cannot occur either.\par

\bigskip
{\bf Case 2c} \quad $\dim D_1 = 3$.\par

\bigskip
If $D_1$ has dimension 3, then $D_1$ is a component of $\loc(W^1_x)$ by proposition \ref{iowifam};
therefore $\dim \loc(W^2)_{\loc(W^1_x)} \ge 4$ by lemma \ref{locy}. 
Let $D_2$ be a component of $\loc(W^2)_{\loc(W^1_x)}$; since $W^2$ is not covering, $D_2$ is
a divisor in $X$ and, by lemma \ref{numeq}, $N_1(D_2)=<[W^1],[W^2]>$.\par

\medskip
Suppose that $D_2 \cdot V > 0$; then up to exchange $W^3$ and $W^4$ we have $D_2 \cdot W^3 >0$, 
hence $\loc(W^3)_{D_2}$ is nonempty and, by lemma \ref{locy}, $\dim \loc(W^3)_{D_2}= 5$, a contradiction.\par

\medskip
So we have $D_2 \cdot V = 0$; in this case $D_2$ must be negative on one of the
$W^i$, but, since $N_1(D_2)= <[W^1],[W^2]>$,  $[W^i]$ must belong to $\Pi$.\\
In particular for every pair $(W^i,W^{i+1}) \in \B_{\Pi}$ we have that $(D_2 \cdot W^i)(D_2 \cdot W^{i+1}) <0$.\\
Moreover, if $D_2 \cdot W^i <0$, then $D_2= \loc(W^i)$; in fact, if $\dim \loc(W^i)=3$ then 
we can apply  lemma \ref{locy} (a) and get $\dim \loc(W^2)_{\loc(W^i_x)}=5$, a contradiction.\par

\medskip
Let $T$ be the union of $\loc(W^i, W^{i+1})$ for $(W^i,W^{i+1}) \in \B_\Pi$, let $z \in T$ and
let $x$ be a point outside $T$; since $X$ is rc$\V$ connected
we can join $x$ to $z$ with a chain of cycles in $\V$;
let $\Gamma$ be the first
irreducible curve in the chain which meets $T$.\par
\medskip
First of all we note that $\Gamma$ cannot meet $D_2$; in fact, since $D_2\cdot V= 0$, $\Gamma$
would be a curve in a family $W^i$ whose class does not lie in the plane $\Pi$,
so that, by lemma \ref{locy} $\dim \loc(W^i)_{D_2} =5$, a contradiction.\par
\medskip
Therefore $\Gamma$ meets a component $D_i \not =D_2$ of the locus of a family $W^i$ of a pair in $\B_\Pi$
such that $\loc(W^{i+1})= D_2$ and such that $D_2 \cap D_i \not = \emptyset$.\\
If $D_i$ has dimension four, then we go back to case 2b, so we can assume that $\dim D_i=3$,
i.e. without loss of generality that $D_i=D_1$.\par
\medskip
By construction, $\Gamma$ cannot belong to a family $W^i$ whose class is contained in the plane $\Pi$
and is not proportional to $V$; on the other hand, if $\Gamma$ belongs to a family $W^i$
whose class is not contained in $\Pi$, then, by lemma \ref{locy} $\dim \loc(W^i,W^{i+1})_{D_1} = 5$, 
a contradiction.\par

\medskip
It follows that either $\Gamma$ belongs to an unsplit family $\alpha V$ whose
numerical class is proportional to $V$ or $\Gamma$ belongs to $V$.\par

\medskip
In the first case $\loc(\alpha V)_{D_1}$ is a divisor $D'$ such that $N_1(D')=<[W^1],[\alpha V]>$;
if  $D' \cdot V > 0$, then we can assume that $\loc(W^3)_{D'}$ is nonempty and so \linebreak
$\dim \loc(W^3)_{D'} =5$, a contradiction.\\
Therefore $D' \cdot V = 0$, but, since $D'$ meets $D_1$ and $D' \not \supset D_1$ then $D' \cdot W^1 > 0$
and $D' \cdot W^2 < 0$, so $D'=D_2$ and the curve is contained in $T$, a contradiction.\par

\medskip
Finally, if $\Gamma$ belongs to $V$, we use the following 
\begin{lemma}\label{nenè}
Let $C$ be an irreducible curve in $V$.
Then either $C \subset \loc(V_x)$ for some $x$ such that $V_x$ is unsplit or $C \subset \loc(W^i)$
for some unsplit family $W^i$ such that $[V] = [W^i] + [W^{i+1}]$.
\end{lemma}

\begin{proof}
If there exists a point $x \in C$ such that $V_x$ is unsplit, then we are in the first case. Otherwise,
for every $x \in C$ there passes a reducible cycle $C^i_x +C^j_x \in \V$. Since the families such that
$[W^i]+[W^{i+1}]= [V]$ are only a finite number, it follows that $C \subset \loc(W^i)$ for some $i$. \end{proof}

\medskip
We thus have two possibilities for $\Gamma$: either $\Gamma \subset \loc(V_x)$, with $V_x$ unsplit, so
$\loc(V_x) \cap D_1 \neq \emptyset$ and therefore $\dim \loc(V_x) \cap D_1 \ge 1$, a contradiction
because $N_1(\loc(V_x))=<[V]>$ and $N_1(D_1)=<[W^1]>$,
or $\Gamma \subset \loc(W^i)$ with $[W^i] \not \in \Pi$; in this case $\loc(W^i,W^{i+1})_{D_1}$
is nonempty and by lemma \ref{locy} $\dim \loc(W^i,W^{i+1})_{D_1} = 5$, 
a contradiction. 

\bibliographystyle{amsplain}

\end{document}